\documentclass[11pt]{article}

\usepackage[latin1]{inputenc}
\usepackage[english]{babel}
\usepackage{amsmath,amssymb, enumerate,array,amsthm}
\usepackage[dvips]{epsfig}
\usepackage{amsfonts}
\usepackage{latexsym}
\usepackage{mathrsfs}

\begin {document}

\title{A characterization of balanced episturmian sequences} 
\author{Genevi\`eve Paquin\thanks{with the support of NSERC (Canada)}   \thanks{Laboratoire de combinatoire et d'informatique math\'ematique, Universit\'e du Qu\'ebec \`a Montr\'eal, CP. 8888 Succ. Centre-Ville, Montr\'eal, (QC) CANADA, H3C 3P8, paquin@lacim.uqam.ca }
 \and  Laurent Vuillon\thanks{Laboratoire de math\'ematiques, CNRS UMR 5127, Universit\'e de Savoie, 73376 Le Bourget-du-lac cedex, France, Laurent.Vuillon@univ-savoie.fr}
 \\
}

\date{\today}

\maketitle

\sloppy

\newtheorem{thm}{Theorem}[section]
\newtheorem{prop}[thm]{Proposition}
\newtheorem{cor}[thm]{Corollary}
\newtheorem{lem}[thm]{Lemma}

\newtheorem{exa}[thm]{Example}
\newtheorem{exas}[thm]{Examples}
\newtheorem{nota}[thm]{Notation}
\newtheorem{defi}[thm]{Definition}
\newtheorem{rem}[thm]{Remark}

\newcommand{\QED}{\rule{1ex}{1ex} \par\medskip}

\newcommand{\nl}{\par\medskip\noindent}
\newcommand{\Proof}{{\nl\it Proof.\ }}

\newcommand{\N}{\mathbb N}
\newcommand{\Z}{\mathbb Z}
\newcommand{\A}{\mathcal{A}}
\newcommand{\Pal}{\textnormal{Pal}}
\newcommand{\Ult}{\textnormal{Ult}}
\newcommand{\bprop}{\begin{prop}}
\newcommand{\eprop}{\end{prop}}
\newcommand{\bcor}{\begin{cor}}
\newcommand{\ecor}{\end{cor}}
\newcommand{\blem}{\begin{lem}}
\newcommand{\elem}{\end{lem}}

\begin{abstract}  It is well known that Sturmian sequences are the aperiodic sequences that are balanced  over a 2-letter alphabet. They are also characterized by their complexity: they have exactly $(n+1)$ factors of length $n$. One possible generalization of Sturmian sequences is the set of infinite sequences over a $k$-letter alphabet, $k \geq 3$, which are closed under reversal and have at most one right special factor for each length. This is the set of episturmian sequences. These are not necessarily balanced over a $k$-letter alphabet, nor are they  necessarily aperiodic. In this paper, we characterize balanced episturmian sequences, periodic or not,  and prove Fraenkel's conjecture for the class of episturmian sequences. This conjecture was first introduced in number theory and has remained unsolved for more than 30 years. It states that for a fixed $k> 2$, there is only one way to cover $\Z$ by $k$ Beatty sequences. The problem can be translated to combinatorics on words: for a $k$-letter alphabet, there exists only one balanced sequence up to letter permutation that has different letter frequencies. 
\end{abstract}

\section{Introduction}

Sturmian sequences are exactly the aperiodic balanced sequences over a 2-letter alphabet  \cite{bs,mh}. A sequence $s$ is balanced if for every letter $a$, the number of $a$'s in any two $n$-length factors  differ by at most 1, for any $n$. Sturmian sequences are also characterized by their number of $n$-length factors: they always have $(n+1)$ factors of length $n$, for every $n$. For sturmian sequences, the two conditions are equivalent. There are two different generalization of Sturmian sequences for alphabets of cardinality $k \geq  3$. The first one is the set of balanced sequences \cite{jb,rt1,lv} and the second one is the set of sequences which are closed under reversal and have at most one right special factor for each length. These are called episturmian sequences, and have been extensively studied by Justin and Pirillo \cite{djp,jj,jv}. It is interesting to note that the two notions coincide for Sturmian sequences, which are both aperiodic episturmian  and aperiodic balanced \cite{bs} sequences over a 2-letter alphabet.  Nevertheless, when the alphabet has 3 letters or more, the two notions no longer coincide.  In particular, episturmian sequences are generally unbalanced over a $k$-letter alphabet, for $k \geq 3$. Thus, a natural question is to characterize which sequences are both episturmian and balanced. We pursue this question here.

We show that there are exactly three different kinds of balanced episturmian sequences, and among them, only one has different letters frequencies. Moreover, this characterization gives a partial proof of the Fraenkel's conjecture  \cite{af1,af2,rt1,rt2}. He has conjectured that for a fixed $k >2$, there is only one covering of $\mathbb{Z}$ by $k$ sequences of the form $\lfloor an+b \rfloor$.  These are called the Beatty sequences. In combinatorics on words, the conjecture can be stated as follows: there is only one  balanced sequence over a $k$-letter alphabet with different letter frequencies, up to letter permutation. In particular, we prove Fraenkel's conjecture in this paper for the class of episturmian sequences. Fraenkel's sequence $(Fr_3)^{\omega}$ for a 3-letter alphabet is $Fr_3=1213121$ and Fraenkel's sequence $(Fr_k)^{\omega}$ for a $k$-letter alphabet $\A =\{1,2,\dots,k\}$ is $Fr_k=Fr_{k-1} k Fr_{k-1} .$ The conjecture is verified for $k=3,4,5,6$  according to the work of  Altman, Gaujal, Hordijk and Tijdeman \cite{agh,rg,ph}. The case $k=7$ has been recently settled by Bar\'at and Varj\'u \cite{bv}.
Many cases have been proved by Simpson \cite{rs}. Graham and O'Bryant \cite{go} have generalized the conjecture to exact $k$-fold coverings and they have proved special cases of the generalized conjecture. 

One interesting corollary of our main result is that Arnoux-Rauzy sequences \cite{ar} are never  balanced episturmian sequences, since every balanced sequence is ultimately periodic.

In this paper, we first recall basic definitions and notation for combinatorics on words, as well as some useful results about episturmian sequences. Then, we show that balanced standard episturmian sequences are described by one of the three following directive sequences:
\begin{itemize}
\item [a)] ${\displaystyle \Delta(s)= 1^n23\dots (k-1)k^\omega;} $ 
\item [b)]   ${\displaystyle \Delta(s)=12\dots (k-1) 1 k \dots (k+\ell-1) ({k+\ell })^\omega}$;
\item [c)] ${\displaystyle \Delta(s) = 123\dots k 1^\omega .}$\end{itemize}
As a result, since episturmian sequences have the same language as standard episturmian sequences, we prove a similar characterization for non standard episturmian sequences. Finally, considering letters  frequencies in the balanced episturmian sequences, we prove Fraenkel's conjecture for the class of episturmian sequences. 

\section{Preliminaries}

Let $\A$ denote a finite alphabet. A {\it finite word} $w$ is an element of the free mono\" \i d $\A ^*$. The $i$-th letter of $w$ is denoted  $w_i$. If $w=w_1w_2\dots w_n$, with $w_i \in \A$, the length of $w$ is $n$ and we write $|w|=n$. By convention, the empty word is denoted $\varepsilon$ and its length is 0. We define the set of non empty finite words as $\A^+=\A^* \setminus \{\varepsilon\}$, and $\A ^\omega$ denotes the set of right infinite words over the alphabet $\A$, also called {\it sequences} for short. Then $\A ^\infty=\A^* \cup \A^\omega$ is the set of finite and right infinite words. 

A word $w \in \A^\infty$ is {\it balanced} if for all factors $u$ and $v$ of $w$ having 
the same length, one has $||u|_a-|v|_a|\leq 1$ for every $a \in \A$.  A word $w \in \A^\infty$ is {\it ultimately periodic} of period $n\in \N$ if $w_i=w_{i+n}$ $\, \, \forall i  \geq \ell$ and $\ell \in \N$. If $\ell=1$, then $w$ is {\it purely periodic}. 

The number of occurences of the letter $a \in \A$ in $w$ is denoted $|w|_{a}$. For a finite word $w$, the frequency of the letter $a$ is defined by $f_a(w)=|w|_a/|w|$. Note that we will compute the frequencies only for ultimately periodic and balanced sequences. Thus, the frequencies always exist (see \cite{agh}). The {\it reversal} of the finite word $w=w_1w_2\dots w_n$ is $\widetilde{w}=w_nw_{n-1} \dots w_1$ and  if $\widetilde{w}=w$, then $w$ is said to be a {\it palindrome}. A finite word $f$ is a factor of $w \in \A^\infty$ if $w=pfs$ for some $p \in \A^*, s \in \A^\infty$. If $p=\varepsilon$ (resp. $s=\varepsilon$), $f$ is called a {\it prefix} (resp. a {\it suffix}) of $w$. Let $u=as$, $a \in \A$ and $s \in \A^\infty$, then, $a^{-1}u=s$. The {\it palindromic right closure} of $w \in \A^*$ is the shortest palindrome $u=w^{(+)}$ with $w$ as prefix. 

The set of factors of $s \in \A^\omega$ is denoted $F(s)$ and $F_n(s)=F(s)\cap \A ^n$ is the set of all factors of $s$ of length $n \in \N$. A factor $f$ of $s$ is {\it right} (resp. {\it left}) {\it special} in $s$ if there exist $a,b \in \A$, $a \neq b$, such that $fa, fb \in F(s)$ (resp. $af, bf \in F(s)$). The {\it alphabet} of $s$ is $\textnormal{Alph}(s)=F(s) \cap \A$ and $\Ult(s)$ is the set of letters occuring infinitely often in $s$. The complexity function is given by $p(n)=F_n(s)$ and is the number of factors of $s$ of length $n \in \N$. Two words $w_1,w_2 \in \A^\infty $ are said equivalent if they have the same set of factors: $F(w_1)=F(w_2)$.

A sequence $s \in \A^\omega$ is {\it episturmian} if $F(s)$ is closed under reversal and $s$ has at most one right special factor of length $n$ for each $n \in \N$. Moreover, an episturmian sequence is {\it standard} if all its left special factors are prefixes of it.  In \cite{djp}, the authors define episturmian sequences using standard episturmian ones as follows:
\begin{defi} An infinite word $t$ is episturmian if $F(t)=F(s)$ for some standard episturmian sequence $s$.   
 \end{defi}

These are characterized by the following definition of standard episturmian sequences given in \cite{djp}:
\begin{defi} \label{standepis} A sequence $s$ is standard episturmian if there exists an infinite sequence $\Delta=x_1 x_2 ....$  , $x_i \in \A$, such that each of the words $u_n$ defined by  $u_{n+1}=(u_nx_n)^{(+)}$, $n \geq 1$, with $u_1=\varepsilon$,  is a prefix of $s$.


\end{defi} 

In Definition \ref{standepis}, the word $\Delta(s)$ is called the {\it directive sequence} of the standard episturmian sequence $s$ and we write $s=\Pal(\Delta(s))$. 

\begin{nota} \textnormal{(\cite{jj})} Let $w=x_1x_2\dots x_n$, $x_i \in \A$, and $u_1=\varepsilon$, \dots, $u_{n+1}=(u_nx^n)^{(+)}$, be the palindromic prefixes of $u_{n+1}$. Then \textnormal{Pal}$(w)$ denotes the word $u_{n+1}$. 
\end{nota}

Recall from \cite {jj} a useful property of the operator $\Pal$. It will be used in almost all our proofs in the next section:
\blem \label{lemjj} Let $x \in \A$. If $w$ is $x$-free, then $\Pal(wx)=\Pal(w)x\Pal(w)$. If $x$ occurs in $w$ write $w=w_1xw_2$ with $w_2$ $x$-free. Then, the longest palindromic prefix of $\Pal(w)$ which is followed by $x$ in $\Pal(w)$ is $\Pal(w_1)$ whence easily $\Pal(wx)=\Pal(w)\Pal(w_1)^{-1}\Pal(w)$. 
\elem

{\exa Let $w=\Pal(123)=1213121$. Then, $\Pal(123\cdot 4)=\Pal(123) \cdot 4\cdot \Pal(123) = 121312141213121$ and $\Pal(123\cdot 2)=\Pal(123)\cdot \Pal ^{-1}(w_1)\cdot \Pal(123)=1213121 (1)^{-1}1213121=1213121213121$, with $w_1=1$ and $w_2 =3$. }

The directive sequence allows to construct easily standard episturmian sequences: 

{\exa Over the alphabet $\A=\{1,2,3\}$, the Tribonacci sequence $t$ \textnormal{(}see \textnormal{\cite{ar}}\textnormal{)}, a standard episturmian sequence, has the directive sequence $\Delta(t)=(123)^\omega$ and then, $u_1=\varepsilon$, $u_2=\underline 1$, $u_3=(1\underline 2)^{(+)}=1\underline 21$, $u_4=(121\underline 3)^{(+)}=121\underline 3121$, $\dots$, $t=\underline 1\underline 21\underline 3121\underline 1213121\underline 2131211213121 \dots$
}

{\rem For clarity issue, we underline the letters of the directive sequence in the corresponding episturmian sequence. 
}

A standard episturmian sequence $s \in \A^\omega $ or any equivalent episturmian sequence is said to be $\mathcal{B}$-{\it strict} if $\Ult(s)=\textnormal{Alph}(s)=\mathcal{B} \subseteq \A$: every letter in $\mathcal{B}=\textnormal{Alph}(s)$ occurs infinitely many times in its directive sequence $\Delta(s)$. In particular, the $\mathcal{A}$-strict episturmian sequences correspond to the Arnoux-Rauzy sequences (see \cite{ar}).

\section{Balanced episturmian sequences}

In this section, we give a characterization of the balanced episturmian sequences. We first study balanced standard episturmian sequences and from that characterization, as standard episturmian sequences have the same language as episturmian sequences, we characterize more generally the balanced episturmian sequences. 

Let consider the first repeated letter of the directive sequence $\Delta(s)$ of a balanced standard episturmian sequence $s$. By {\it first repeated letter}, we mean the first letter which occurs twice in the shortest prefix of $\Delta(s)$.

{\exa In $\Delta(s)=12321 \dots$, the first repeated letter is $2$. 
}

{\rem In this paper, in order to show that a word is unbalanced, we will always give two factors $f$ and $f'$ of the same length, with $||f|_{f_1}-|f'|_{f_1}| \geq 2$}: the unbalance is over the first letter of $f$, namely $f_1$.\\

Let start with introductive examples. 

{\exas ${  }$
\begin{itemize} 
\item [1)] Let $s$ be a standard episturmian sequence with the directive sequence $\Delta(s)=1232\dots$. Then, $$s=\Pal (1232\dots)=\underline 1 \underline 21 \underline 3121\underline 2 13121\cdots,$$ which contains the factors $212$ and $131$. Thus, $s$ is unbalanced over the letter $2$.
\item [2)] Let $t$ be a standard episturmian sequence with the directive sequence $\Delta(t)=12131\dots $. Then $$t=\Pal (12131\dots)= \underline 1 \underline 21 \underline 121 \underline 3 121121 \underline 1213121121\dots ,$$  which contains the factors $11211$ and $21312$. Thus, $t$ is unbalanced over the letter $1$. 
\item [3)] Let $u$ be a standard episturmian sequence with the directive sequence $\Delta(u)=12341\dots$. Then $$u=\Pal(12341\dots)=\underline 1\underline 21\underline 3121\underline 41213121\underline 121312141213121\dots,$$ which is a balanced prefix.
\end{itemize}
}

It seems that the satisfaisability  of the balance condition depends on where the repeated letters occur. Proposition \ref{lem1} characterizes directive sequences with the first repeated letter different from the first letter, while Proposition \ref{lem2} characterizes the directive sequences $\Delta(s)=11z$, with $z \in \A^\omega$. 

\blem \label{lemg}Let $\Delta(s)=x\alpha ^\ell y$ be the directive sequence of a balanced standard episturmian sequence $\alpha \in \A$, $x \in \A^+$, $y \in \A^\omega$ and $\ell \geq 2$. If $\alpha$ is the first repeated letter, then $x_i \neq y_j$, $\forall i, j$.  
\elem

\Proof Let suppose there exists $\beta \neq \alpha \in \A $ such that $\Delta(s)=x' \beta x''  \alpha ^\ell y'  \beta y'' $, with $y'_i \neq x_j$, $\forall i, j$. Let $p=\Pal (x' \beta x'')$. There are 3 cases:
\begin{itemize}
\item [a)] If $x' \neq \varepsilon$: then $$s=p(\underline \alpha p)^\ell \dots (p\alpha)^\ell p \underline \beta p_1 \dots,$$ which contains the factors $\alpha p \alpha$ and $p \beta p_1$, $p_1 \neq \alpha$, $p_1 \neq \beta$. Thus, $s$ is unbalanced. 
\item [b)] If $x' = \varepsilon$ and $x'' \neq \varepsilon$: then $$s=(p(\underline \alpha p)^\ell \dots (p\alpha)^\ell p )^2 \dots ,$$ which contains $\alpha p \alpha$ and $pp_1p_2$, $p_1 =\beta \neq \alpha$, $p_2=x_1''\neq \alpha$. Then $s$ is unbalanced. 
\item [c)] If $x'=\varepsilon$ and $x''=\varepsilon$: then $$s=\underline \beta (\underline \alpha \beta)^\ell \dots \beta \underline \gamma \beta \dots,$$ which contains $\alpha \beta \alpha$ and $ \beta \gamma \beta$, $\gamma \in \A$. Since $s$ is over at least a 3-letter alphabet, $\gamma$ exists and at its first occurrence it is preceded and followed by $\beta$. Then, $s$ is unbalanced.  
\end{itemize}
\QED

\bprop \label{lem1} Let $\Delta(s)$ be the directive sequence of a balanced standard episturmian sequence $s$. If the first repeated letter $k$ of the directive sequence is not the first letter of $s$, then the directive sequence can be written as $\Delta(s)=12\dots (k-1)k^\omega$.
\eprop
\Proof Let $\Delta(s)=xky kz$ be the directive sequence of a balanced standard episturmian sequence $s$,  with $x \in \A ^+$, $y \in \A^*$, $z \in \A^\omega$, $|xky|_{\alpha} \leq 1$ , $\forall \alpha \in \A$ (since $k$ is the first repeated letter).  Let $p=\Pal(x)$ and let suppose $y \neq \varepsilon$. Then
$$s=p\underline kp\underline y_1 p kp \dots pkp y_1 pkp \underline kp\dots$$ which contains the factors $kp k$ and $py_1p_1$, $p_1 \neq k$ and $y_1 \neq k$ (since $|xky|_{\alpha} \leq 1$, $\forall \alpha \in \A$). Then $s$ is unbalanced over $k$. Thus, it follows that $y=\varepsilon$. Consequently, $\Delta(s)=x k^2 z$. We rewrite $\Delta(s)=x k^\ell z'$, with $z'_1 \neq k$, and $\ell \geq 2$. Then,
$$s= p(\underline kp)^\ell \underline  z_1'p_1 \dots$$ which contains the factors $kp k$ and $p z_1'p_1$, with $z_1' \neq k$ and $p_1 =x_1 \neq k$. Lemma \ref{lemg} insures that it is the first occurence of the letter $z_1'$ and then, $z_1'$ is necessarily followed by $p_1$. If $z_1'\neq k$ exists, $s$ is unbalanced. Then $z'=k^\omega$ and $\Delta(s)=12 \dots (k-1) k^\omega$. 
\QED

\bprop \label{lem2} Let $\Delta(s)$ be the directive sequence of a balanced standard episturmian sequence $s$. If $\Delta(s)=1^\ell z$, with $z \in \A^\omega$, $z_1 \neq 1$ and $\ell \geq 2$, then $\Delta(s)=1 ^\ell  23\dots (k-1) k^\omega$.
\eprop

\Proof  Let $\Delta(s) = 1 ^\ell z$ be the directive sequence of a balanced standard episturmian sequence, with $z_1 \neq 1$ , $\ell \geq 2$. Let suppose  $|z|_{1} > 0$. Then we get 
 $\Delta(s)=1^\ell z' 1 z''$, with $z' \neq \varepsilon$ and $|z'|_{1} =0$. Then, 
 $$s= \underline 1 ^\ell \underline z'_11 ^\ell \dots 1 ^\ell z'_11 ^\ell \underline 1 z'_1   \dots,$$ which contains the factors $ z'_11 ^{\ell+1}z'_1$ and $1 ^\ell  \alpha 1 ^2$, with $\alpha \in \A$, $\alpha \neq 1$ and $\alpha \neq z'_1$. There is at least one letter $\alpha$ in $z'$ or $z''$ distinct from $z'_1$ and $1$, since $s$ is over at least a 3-letter alphabet, and at its first occurence, $\alpha$ is preceeded and followed by $1^\ell$. Then $s$ is unbalanced. It follows that $|z|_{1}=0$. Since the alphabet is finite, there is at least one letter distinct from $1$ which occurs twice. Let consider the first repeated one, namely $\gamma$. Then, $\Delta(s)=1 ^\ell u \gamma v \gamma w$, with $|u\gamma v|_{\alpha} \leq 1$, $\forall \alpha \in \A$ and $|u\gamma v \gamma w|_{1}=0$. Let suppose $v \neq \varepsilon$ and let $p=\Pal (1 ^\ell u)$. Then, 
 $$s=p\underline \gamma p \underline v_1 p\gamma p \dots p \gamma p \underline \gamma \dots $$
 which contains the factors $\gamma p \gamma$, $pv_1p_1$, $v_1=1\neq \gamma$ and $p_1 \neq \gamma$. We conclude that $v=\varepsilon$. 
Let now consider $\Delta(s)= 1 ^\ell u \gamma ^2 w$, which we rewrite as $\Delta(s)=1 ^\ell u \gamma ^m w'$, $m \geq 2$ and $w_1' \neq \gamma$. Then,
$$s=p(\underline \gamma p)^m \underline w'_1 p_1 \dots$$ which contains the factors $\gamma p \gamma$ and $pw'_1p_1$, $p_1 \neq \gamma$, $w_1' \neq \gamma$. Then, $w'=\gamma ^\omega$. We conclude that $\Delta(s)=1 ^\ell  23  \dots (k-1)\gamma^\omega=1 ^\ell  23  \dots (k-1)k^\omega$.   \\
\QED

Propositions \ref{lem5} and \ref{lem3} characterize the directive sequences $\Delta(s)=1y1z$, with $y \neq \varepsilon$, $|y|_1=0$ and $1$, the first repeated letter. Two technical lemmas are required:

\blem \label{lem2iemedouble} Let $\Delta(s)=1y1z$ be the directive sequence of a balanced standard episturmian sequence, with $1$ the first repeated letter, $|y|_1=0$, $y \in \A^+$, $z \in \A^\omega$. Then, $z_i \neq y_j$, $\forall i, j$. 
\elem

\Proof Let $\Delta(s)=1y'\alpha y'' 1 z' \alpha z''$, with $\alpha \in \A$, $z'_i \neq y_j$, $\forall i, j$ (the second repeated letter is $\alpha$). Let $p=\Pal(1y')$. Then, there are 2 cases:
\begin{itemize}
\item [a)] If $y' \neq \varepsilon$, then $$s=(p\underline \alpha p \dots p\alpha p)^2 \dots (p\alpha p \dots p \alpha p)^2 \underline \alpha \dots ,$$ 
which contains the factors $\alpha p \alpha$ and $pp_1p_2=p1y'_1$ and $1\neq \alpha$, $y_1' \neq \alpha$.  Then $s$ is unbalanced.
\item [b)] If $y'=\varepsilon$, then $s=(\underline 1 \underline \alpha 1 \dots 1 \alpha 1)^2\dots (1\alpha 1 \dots 1 \alpha 1)^2 \underline \alpha \dots $, which contains the factors $\alpha 1\alpha$ and $1\beta 1$, $\beta \in \A \setminus \{1, \alpha \}$ is necessarily preceded and followed by the letter 1. The letter $\beta$ exists since $s$ is over at least a 3-letter alphabet. Then $s$ is unbalanced.
\end{itemize}
It follows that $z_i \neq y_j$, $\forall i,j$. 
\QED

\blem \label{lem4} Let $\Delta(s)= {1} y {1} z$ be the directive sequence of a balanced standard episturmian sequence $s$, with 1 the first repeated letter, $y \in \A^+$, $|y|_{1} =0$, $z_1 \neq 1$ and $z \in \A^\omega$. Then, $|z|_{1} =0$. 
\elem

\Proof Let suppose $|z|_{1} \geq 1$ (there is a third ${1}$ in the directive sequence): $\Delta(s)={1} y {1} z' {1} z''   $, with $|z'|_{1}=0$, $z' \neq \varepsilon $. Let $p=\Pal({1} y)$. Then
$$s=\underline p^2  \underline z'_1 p^2 \dots p^2 z'_1 p^2 \underline p z'_1 p^2 \dots  ,$$ which contains the factors ${1} p {1} $ (in $p^3$) and ${1} ^{-1}p z_1' p_1 p_2$, where $z_1' \neq p_1={1}$ and $z_1' \neq p_2=y_1$ (insures by Lemma \ref{lem2iemedouble}).

If the directive sequence is $\Delta(s)= {1} y {1} z$, then $z$ does not contain any other ${1}$.
\QED

Lemma \ref{lem4} insures there is no more occurence of $1$ in $\Delta(s)$.  Let see some examples where there is at least a second repeated letter:

{\exa Let $s$ be a standard episturmian sequence with directive sequence $\Delta(s)=12321\dots $. Then, $$s=\underline 1 \underline 21 \underline 3121 \underline 213121 \underline 1213121213121\dots$$ which contains the factors $212$ and $131$. Thus, $s$ is unbalanced over the letter $2$.
}

{\exa Let $t$ be a standard episturmian sequence with directive sequence $\Delta(t)=12312\dots$. Then, $$t=\underline 1 \underline 21 \underline 3121 \underline 1213121 \underline 213121\dots $$ which contains the factors $212$ and $131$. Thus, $t$ is unbalanced over the letter $2$.
}

In the following proposition, we describe how the letters different from $1$ can be repeated: 

\bprop \label{lem5} Let $\Delta(s)={1} y {1} z$ be the directive sequence of a balanced standard episturmian sequence, with $|y|_1=|z|_1=0$, $y_i \neq y_j$, $\forall i \neq j$, $y \neq \varepsilon$. Then,  $\Delta(s)={1}23 \dots k {1} (k+1) \dots (k+\ell) (k+\ell +1)^\omega$. 
\eprop

\Proof Let consider $\alpha \in \A$, the first repeated letter distinct from ${1}$. Then, there are 2 cases to consider:
\begin{itemize}
\item [a)] $\Delta(s)= {1} y' \alpha y'' {1}  z' \alpha z''$, with $|y'|_\alpha=|y''1z'|_\alpha=0$.  Impossible, by Lemma \ref{lem2iemedouble}.
\item [b)] $\Delta(s)={1} y' {1} z' \alpha z'' \alpha z'''$, with $z'' \neq \varepsilon$, $|y'1z'|_\alpha=|z''|_\alpha=|z'\alpha z'' \alpha z''|_{1}=0$. Let $p=\Pal({1} y' {1} z')$. Then
$$s=p\underline \alpha p\underline z_1''p\alpha p\dots p \alpha p z_1''p\alpha p \underline \alpha \dots$$ with factors $\alpha p \alpha$ and $pz''_1p_1$. We  conclude that $z''=\varepsilon$. 
\end{itemize}
The only possibility is $\Delta(s)={1} y' {1} z' \alpha ^2 z'''$, which we rewrite as $\Delta(s)={1} y' {1} z' \alpha ^\ell z''''$, $l \geq 2$, $z''''_1 \neq \alpha$. Let $p=\Pal({1} y' {1} z')$. Then, 
$$s=p(\underline \alpha p)^\ell \underline  z_1''''{p_1}\dots$$ with factors $\alpha p \alpha$ and $p z''''_1{p_1}$, with $z_1'''' \neq \alpha$, $p_1=1\neq \alpha$. It follows that $z''''=\alpha ^\omega$, then $z'''=\alpha ^\omega$ and finally, $\Delta(s)={1} 2 \dots k {1} (k+1) \dots (k+\ell)\alpha^\omega={1} 2 \dots k {1} (k+1) \dots (k+\ell)(k+\ell +1)^\omega$.\QED

\bprop \label{lem3} Let $\Delta(s)=1y {1}^\ell z$ be the directive sequence of a balanced standard episturmian sequence $s$, with $\ell \geq 2$, $z_1 \neq 1$, $y \in \A^+$, $|y|_{1}=0$, $z \in \A^\omega$ and $y_i \neq y_j$, $\forall i \neq j$. Then, $\Delta(s)= 1 2 3 \dots k (1) ^\omega$. 
\eprop 

\Proof Let $p=\Pal(1 y)$. Then,
$$s= p^{\ell +1}\dots p^{\ell +1}\underline z_i1y_1\dots$$ which contains the factors $1 p 1$  (in $p^3$) and $1 ^{-1} p z_i 1 y_1$, with $y_1 \neq 1$, $z_i$ the first letter of $z$ such that $z_i \neq 1$ and $z_i \neq y_1$. Then, $s$ is unbalanced. It follows that $z=1 ^\omega$ and $\Delta(s)=1 2 3 \dots k (1) ^\omega$. \QED

We have now considered every possibility of directive sequences for a balanced standard episturmian sequence. Theorem \ref{3cas} summarizes the previous propositions. 

\begin{thm} \label{3cas} Any balanced standard episturmian sequence $s$ has a directive sequence in one of the three following families of sequences:
\begin{itemize}
\item [a)] ${\displaystyle \Delta(s)=1^n \left ( \prod_{i=2}^{k-1}  i \right ) (k)^\omega  = 1^n23\dots (k-1)(k)^\omega,}$ with $n \geq 1$; 
\item [b)]  ${ \displaystyle \Delta(s)=\left (  \prod_{i=1}^{k-1}i \right ) 1 \left ( \prod_{i=k}^{ k+ \ell-1}i \right ) (k+\ell)^\omega = 12\dots (k-1) 1 k \dots (k+\ell-1) (k+\ell )^\omega}$, with $k \geq 2$, $l \geq 1$.
\item [c)]  ${ \displaystyle \Delta(s) = \left (\prod _{i=1}^k i\right ) (1)^\omega =123\dots k (1)^\omega}$, with $k \geq 3$; 
\end{itemize}
\end{thm}

\Proof Proposition \ref{lem1} implies a) for $n=1$ and Proposition \ref{lem2} implies a) for $n \geq 2$,  while b) (resp. c)) follows from Proposition \ref{lem5} (resp. Proposition \ref{lem3}). \QED




\begin{rem}
Notice that all the directive sequences given in Theorem \ref{3cas} yield 
balanced standard episturmian sequences. The sequences given in $a)$ can 
be written as:
$$ (1^n21^n31^n21^n41^n21^n31^n21^n \dots  
1^n21^n31^n21^n41^n21^n31^n21^nk)^\omega.$$ 
This sequence is balanced over the letter 1, since the projection of $s$ is 
$(1^n\alpha)^\omega $ which is balanced, and as the distance between two 
occurrences of the same letter is always the same, it follows from Hubert 
\cite{ph} that it is a balanced sequence.  That is the letter $\alpha$ is 
periodically replaced by $\Pal(2 3\dots k)$.

The sequences given in $b)$ can be written as $Fr_{k-1}(Fr_{k-1}A)^\omega$, 
with $A$ periodically replaced by $\Pal(k (k+1)\dots (k+\ell))$. Then, if the 
letter $i \in \{k, k+1, \dots, k+\ell\}$ then from Hubert, the sequence is 
balanced over this letter. The sequence is balanced over the the letter $i \in \{1,2,\dots, k-1\}$, since it appears in the Fraenkel word $(Fr_{k-1}A Fr_{k-1})
^\omega$.  

Finally, the sequences given in $c)$ are known to be balanced, from the 
Fraenkel's conjecture. 
\end{rem}

Recall the following result:
\begin{thm} \label{thdjp} \textnormal{(\cite{djp}, Theorem 3)} A standard episturmian sequence $s$ is ultimately periodic if and only if its directive sequence $\Delta (s)$ has the form $w\alpha ^\omega$, $w \in \A^*$, $\alpha \in \A$. 
\end{thm}

Then, a standard episturmian sequence can not be both periodic and $\mathcal{A}$-{\it strict}. The next result follows:

\bcor \label{pasAR}Every balanced standard episturmian sequence is ultimately periodic. In other words, none of the Arnoux-Rauzy sequences ($\mathcal{A}$-strict episturmian sequences) are balanced. 
\ecor

\Proof The first part of the corollary is a direct consequence of Theorem \ref{3cas}. The second part follows from the property that an $\A$-strict episturmian sequence can not be ultimately periodic. \QED

\bcor \label{3mots}Any balanced standard episturmian sequence $s$ is in one of the following families:
\begin{itemize}
\item [a)] $s= p (k-1) p \left (kp (k-1) p\right )^\omega  $, with $p=\Pal(1^n2 \dots (k-2))$;
\item [b)] $s=p(k+\ell-1)p\left ((k+\ell)p(k+\ell-1)p\right )^\omega $, with $p=\Pal(123\dots (k-1)1k\dots (k+\ell-2))$;
\item [c)] $s= \left [\Pal(123\dots k)\right ] ^\omega $.
\end{itemize}
\ecor

\Proof It follows from the computation of $\Pal (\Delta(s))$, with $\Delta(s)$, the directive sequences of Theorem \ref{3cas}.\QED

\bprop \label{propfreq} Every balanced standard episturmian sequence $s$ with different frequencies for every letter can be written as in Corollary \ref{3mots} c).
\eprop

\Proof In Corollary \ref{3mots} a), the letters $k$ and $(k-1)$ appear once in the period. Thus, the frequencies of $k$ and $(k-1)$ in $s$ are equal. In b), the same argument over the letters $(k+\ell)$ and $(k+\ell-1)$ holds. In c),  by direct inspection, we find that the period of $s$ has length $|\Pal(123\dots k)|=2^k-1$, and the frequency of the letter $i$ is $2^{k-i}/(2^k-1)$. Thus, the frequencies of two different letters are different. \QED

As every episturmian sequence $t$ has the same language as a standard episturmian sequence $s$, the result of Proposition \ref{propfreq} can be extended to any balanced episturmian sequences. Then, we get the following general result, which is the proof of Fraenkel's conjecture for episturmian sequences:

\begin{thm} Let $s$ be a balanced episturmian sequence over a finite $k$-letter alphabet, $k \geq 3$ with $f_i(s) \neq f_j(s)$, $\forall i \neq j$. In other words, every letter frequency is different. Then, $s= \left [\Pal(123\dots k)\right ] ^\omega$, $1, 2, \dots k \in \A$.   
\end{thm}

{\exa For $k=3, 4,5$, we obtain respectively $s=(1213121)^\omega$, $t=(121312141213121)^\omega$ and $u=(1213121412131215121312141213121)^\omega$.
}

\section{Concluding remarks}

In order to extend our work on balanced sequences, we could investigate two directions.

The first one is in relation to billiard sequences.
In \cite{lv}, the second author surveys balanced sequences and proves that billiard sequences over a $k$-letter alphabet are $(k-1)$-balanced. That is, let $X$ be a billiard sequence over a $k$-letter alphabet then
$$\forall i \in \mathcal{A}, \forall n \in \N, \forall w,w' \in F_n(X) \mbox{ we have } ||w|_i-|w'|_i| \leq k-1.$$
Furthermore, we notice that Fraenkel's sequences are  periodic billiard sequences. Therefore,  in the spirit of our paper, it would be interesting to study the billiard sequences \cite{amst,br} which are balanced.  This class contains at least Fraenkel's sequences, and perhaps some other interesting sequences.

The second direction is directly related to the original Fraenkel's conjecture. In order to prove this conjecture, it will be useful to have the property that balanced sequences over an alphabet with more than 2 letters, with pairwise distinct frequencies of letters, are given by directive sequences. Combining our results with this conjecture would give a proof of Fraenkel's conjecture.

\end{document}